\input amstex

\NoBlackBoxes

\documentstyle{amsppt}

\topmatter

\title{One dimensional topological Galois theory}\endtitle

\author Askold Khovanskii\endauthor

\affil{Department of Mathematics, University of Toronto, Toronto, Canada}\endaffil


\thanks{This work was partially supported by the Canadian Grant No. 156833-17.}\endthanks

\abstract{In the preprint we present an outline of the one dimensional version of topological Galois theory. The theory   studies topological obstruction to solvability  of equations ``in finite terms" (i.e. to their solvabilty by radicals, by elementary functions, by quadratures and so on).
The preprint is based on the author's book on topological Galois theory. It contains definitions, statements of results and comments to them. Basically no proofs are presented.

 The preprint was  written as a part of the comments to a new edition (in preparation) of the classical book ``Integration in finite terms'' by J.F.~Ritt.  }\endabstract



\endtopmatter
\document

\subhead {1. Introduction}\endsubhead  As was discovered by  Camille Jordan  the monodromy group of an algebraic function is isomorphic to the
Galois group of the associated extension of the field of
rational functions. Therefore  the monodromy group is
responsible for the representability of an algebraic function
by radicals (see~[1]).

However, not only algebraic functions have the
monodromy group.  It is defined for any solution of a linear differential equation whose  coefficients are rational  functions and for many more functions, for which the Galois group does not make sense.  It is thus natural to try using the monodromy group
for these functions instead of the Galois group to prove that
they do not belong to a certain Liouvillian class. This particular
approach is implemented in   topological
Galois theory (see~[2]), which has a one-dimensional version and a multidimensional version.

In the one-dimensional version we consider functions from Liouvillian classes  as multi-valued analytic functions of one complex variable.
It turns out that there exist topological restrictions on the way the
Riemann surface of a function from a certain Liouvillian class   can be positioned over the complex plane. If a function does not satisfy these restrictions, then it cannot  belong to the corresponding Liouvillian class.

Besides a geometric appeal, this approach has the following
advantage. Topological obstructions relate to 
branching. It turns out that if a function does not belong to a certain Liouvillian class by topological reasons then it automatically does not belong to a much wider class of functions. 
This wider class can be obtained if we
add to the Liouvillian  class all single valued functions having at most countable set of singularities   and allow them to enter all formulas.

The composition of functions is not
an algebraic operation. In differential algebra, this operation
is replaced with a differential equation describing it.
However, for example, the Euler $\Gamma$-function does not satisfy any
algebraic differential equation. Hence it is pointless to look
for an equation satisfied by, say, the function $\Gamma(\exp\,
x)$  and one can not describe it algebraically  (but the function $y=\exp (\Gamma(x))$ satisfies the equation $y'=\Gamma' y$ over a differential field containing $\Gamma$ and it makes sense in the differential algebra).  The only known results  on non-representability of
functions by quadratures and, say, the Euler $\Gamma$-functions
are obtained by our method.

On the other hand, our method cannot be used to prove that a
particular single valued meromorphic function does  not belong to a certain Liouvillian class.

There are the following topological obstructions to representability of functions by 
generalized quadratures, $k$-quadratures and quadratures.

Firstly, the functions representable by generalized quadratures and, in particular, the
functions representable by $k$-quadratures and quadratures  may have no more than
countably many singular points in the complex plane (see section 6).

Secondly, the monodromy group of a function representable by quadratures is
necessarily solvable (see section 8). There are similar restrictions for for a function representable by generalized quadratures and $k$-quadratures.
However, these restrictions are more involved.
To state them, the monodromy group should be regarded not as an abstract group
but rather as a transitive subgroup in the permutation group.
In other terms, these restrictions make use not only of the monodromy group
but rather of the {\it monodromy pair} of the function consisting of
the monodromy group and the stabilizer of some germ of the function
(see section 9).
\smallskip

One can prove that the only
reasons for unsolvability in finite terms of  Fuchsian linear
differential equations are topological (see section 14). In other words, if there are no topological
obstructions to solvability of a Fuchsian equation by  generalized quadratures (by $k$-quadratures, by quadratures), then this equation is solvable by  generalized quadratures (by $k$-quadratures or  by quadratures respectively).
The proof  is based on a linear-algebraic part of differential Galois theory (dealing with linear algebraic groups and their differential invariants),

\subhead{2. Solvability of equations in finite terms}
\endsubhead
An equation is solvable ``in finite terms'' (or is solvable ``explicitly'') if its solutions belong to a certain class of functions. Different classes of functions correspond to different notions of  solvability in finite terms.

A class of functions can be introduced by specifying a list of {\it basic functions}
and a list of {\it admissible operations}.
Given the two lists, the class of functions is defined as the set
of all functions that can be obtained
from the basic functions by repeated application of admissible operations.
Below, we define Liouvillian classes of functions  in exactly this way.

Classes of functions, which appear in the problems of integrability in finite terms,
contain multivalued functions.
Thus the basic terminology should be made clear.
We work with multivalued functions ``globally'', which leads
to a more general understanding of classes of functions defined by
lists of basic functions and of admissible operations.
A multivalued function is regarded as a single entity.
{\it Operations on multivalued functions} can be defined.
The result of such an operation is a set of multivalued functions;
every element of this set is called a function obtained from the given functions
by the given operation.
A  {\it class of functions} is defined as the set of all (multivalued) functions
that can be obtained from the basic functions by repeated
application of admissible operations.

\subhead{3. Operations on multivalued functions}
\endsubhead
Let us define, for example, the sum of two multivalued functions on a connected Riemann surface $U$.

\definition{Definition 7}
Take an arbitrary point $a$ in $U$, any germ
$f_a$ of an analytic function $f$ at the point $a$ and any germ $g_a$ of an analytic function $g$
at the same point $a$.
We say that the multivalued function $\varphi$ on $U$ generated by the germ $\varphi_a=f_a+g_a$
{\it is representable as the sum of the functions} $f$ and $g$.
\enddefinition

For example, it is easy to see that exactly two functions of one variable
are representable in the form
$\sqrt{x}+\sqrt{x}$, namely, $f_1=2\sqrt{x}$ and $f_2\equiv 0$.
Other operations on multivalued functions are defined in exactly the same way.
{\it For a class of multivalued functions, being stable under addition means that,
together with any pair
of its functions, this class contains all functions representable as their sum.}
The same applies to all other operations on multivalued functions understood
in the same sense as above.

In the definition given above, not only the operation of addition plays a key role but
also the operation of analytic continuation hidden in the notion
of multivalued function.
Indeed, consider the following example.
Let $f_1$ be an analytic function defined on an open subset $V$ of the complex line
$\Bbb C^1$ and admitting no analytic continuation outside of $V$,
and let $f_2$ be an analytic function on
$V$ given by the formula $f_2=-f_1$.
According to our definition, the zero function is representable in the form
$f_1+f_2$ {\it on the entire complex line}.
By the commonly accepted viewpoint, the equality $f_1+f_2=0$ holds inside
the region $V$ but not outside.

Working with multivalued functions globally, we do not insist on the existence of
{\it a common region}, were all necessary operations would be performed on
single-valued branches of multivalued functions.
A first operation can be performed in a first region,
then a second operation can be performed in a second, different region
on analytic continuations of functions obtained on the first step.
In essence, this more general understanding of operations is equivalent to including
analytic continuation to the list of admissible operations on the analytic germs.

\subhead {4. Liouvillian classes of single variable functions}\endsubhead
In this section, we define Liouvillian classes of single variable functions
(for many variables, the corresponding definitions can be found in [ 2]).
We will describe these classes by lists of basic functions and admissible operations.
\smallskip
We will need the list of {\it basic elementary functions}.
In essence, this list contains functions that are studied in high-school and which
are frequently used in pocket calculators.

\subhead {List of basic elementary functions}\endsubhead

1. All complex constants and an independent variable $x$.

2. The exponential, the logarithm and the power $x^\alpha$,
where $\alpha$ is any complex constant.

3. Trigonometric functions:  sine, cosine, tangent, cotangent.

4. Inverse trigonometric functions: arcsine, arccosine, arctangent, arccotangent.
\bigskip
\proclaim{Lemma 1}
Basic elementary functions can be expressed through the exponentials and
the logarithms with the help of complex constants, arithmetic operations
and compositions.
\endproclaim
Lemma 1 can be considered as a  simple exercise. Its proof can be found in [3].

Let us now proceed with the list of {\it classical operations} on functions.

\proclaim {List of  classical  operations}

1. The operation of composition takes functions $f$,$g$ to the function $f\circ g$.

2.  The arithmetic operations take functions $f$, $g$ to the functions $f+g$, $f-g$, $fg$, and $f/g$.

3. The operation of differentiation takes function $f$ to the function $f'$.

4. The meromorphic operation takes functions $f_1,\dots,f_n$ to the function $F(f_1,\dots,f_n)$ where $F$ is a fixed meromorphic function of $n$ complex variables.

5. The operation of integration takes function $f$ to a solution of equation $y'=f$ (the function $y$ is defined up to an additive constant).

6.  The operation of solving algebraic equations takes functions $f_1,\dots,f_n$ to the function $y$ such that $y^n+f_1y^{n-1}+\dots+f_n=0$ (the function $y$ is not quite uniquely determined by functions $f_1,\dots,f_n$ since an algebraic equation of degree $n$ can have $n$ solutions).

7. The operation of solving linear differential equations takes functions $f_1,\dots,f_n$ to the function $y$ such that $y^{(n)}+f_1y^{(n-1)}+\dots+f_n=0$ (the function $y$ is not  uniquely determined by functions $f_1,\dots,f_n$ since an differential equation of order   $n$ has an $n$ dimensional space of solutions).

\endproclaim

We can now return to the definition of Liouvillian classes of single variable functions.

\subhead {Functions  representable by radicals}\endsubhead
List of basic functions: all complex constants, an independent variable $x$.
List of admissible operations:
arithmetic operations and the operation of taking the $n$-th root
$f^{\frac {1}{n}}$, $n=2,3,\dots$, of a given function~$f$.

\subhead {Functions  representable by $k$-radicals}\endsubhead
List of basic functions: all complex constants, an independent variable $x$.
List of admissible operations:
arithmetic operations and the operation of taking the $n$-th root
$f^{\frac {1}{n}}$, $n=2,3,\dots$, of a given function $f$, the operation of solving algebraic equations of degree $\leq k$.
\bigskip

\subhead{Elementary functions}\endsubhead
List of basic functions: basic elementary functions.
List of admissible operations: compositions, arithmetic operations, differentiation.

\subhead {Generalized elementary functions}\endsubhead
This class of functions is defined in the same way as the class of elementary functions.
We only need to add the operation of solving algebraic equations to the list
of admissible operations.

\subhead {Functions  representable by quadratures}\endsubhead
List of basic functions: basic elementary functions.
List of admissible operations: compositions, arithmetic operations, differentiation, integration.

\subhead{Functions  representable by $k$-quadratures}\endsubhead
This class of functions is defined in the same way as the class of functions representable by quadratures.
We only need to add the operation of solving algebraic equations of degree
at most $k$ to the list of admissible operations.

\subhead{Functions  representable by generalized quadratures}\endsubhead
This class of functions is defined in the same way as
the class of functions representable by quadratures.
We only need to add the operation of solving algebraic equations to the
list of admissible operations.

\subhead {5. Simple formulas with complicated topology}
\endsubhead
Developing topological Galois theory I followed the following plan:
\medskip

I. To find a wide class of multivalued functions such that: 
\smallskip

a) it is closed under all classical operations; 

b) it contains all entire functions and all functions from each Liouvillian class; 

c) for functions from the class the monodromy group  is well defined.
\medskip

II. To use the monodromy group instead of the Galois group inside the class.

\medskip

Let us discuss some difficulties that one need to overcome on this
way.

\example{Example} Consider an elementary function $f$ defined by the
following formula:
$$
f(z)=\ln \sum\limits_{j=1}^n\lambda_j\ln (z-a_j)
$$
where $a_{j}$  are different points in the
complex line, and $\lambda_j\in \Bbb C$ are constants.  
\endexample

Let $\Lambda$ denote the additive subgroup of
complex numbers generated by the constants $\lambda_1$,
$\dots$, $\lambda_n$. It is clear that if $n>2$, then for
almost every collection of constants $\lambda_1,\dots,
\lambda_n$, the group $\Lambda$ is everywhere dense in the
complex line.

\proclaim {Lemma 2}
  If the group $\Lambda$ is dense in the complex line, then the elementary function $f$
has a dense set of logarithmic ramification points.
\endproclaim

\demo{Proof}
  Let $g$ be the multivalued function defined by the formula
$$g(z)=\sum\limits_{j=1}^n\lambda_j\ln (z-a_j).$$ Take  a point 
$a\neq a_j$, $j=1,\dots, n$ and let $g_a$ be one of the germs of $g$ at $a$. A loop around the points
$a_1,\dots, a_n$ adds the number $2\pi i\lambda$ to the germ
$g_a$, where $\lambda$ is an element of the group $\Lambda$.
Conversely, every germ $g_a+2\pi i\lambda$, where $\lambda \in
\Lambda$, can be obtained from the germ $g_a$ by the analytic
continuation along some loop. Let $U$ be a small neighborhood
of the point $a$, such that the germ $g_a$ has  a singe-valued analytic continuation $G$  on $U$.  The image $V$
of the domain $U$ under the map $G: U\to\Bbb C$ is open.
Therefore, in the domain $V$, there is a point of the form
$2\pi i \lambda$, where $\lambda \in \Lambda$. The function
$G-2\pi i \lambda$ is one of the branches of the function $g$
over the domain $U$, and the zero set of this branch in the
domain $U$ is nonempty. Hence, one of the branches of the
function $f=\ln g$ has a logarithmic ramification point in $U$.
\enddemo

The set $\Sigma$ of singular points of  the function $f$ is  a {\it countable set} (see   section 6). Under assumptions of Lemma 2 the set $\Sigma$ is  everywhere dense.

It is not hard to verify that  the monodromy group (see section 7)  of the function $f$ has the cardinality of the continuum.  This is not surprising: the
fundamental group $\pi_1(\Bbb C \setminus \Sigma)$ has obviously the
cardinality of the continuum provided that $\Sigma$ is a countable
dense set.

One can also prove that the image of the fundamental group $\pi
_1(\Bbb C\setminus \{\Sigma \cup b\})$ of the complement of the set
$\Sigma\cup b$, where $b\not\in \Sigma$, in the permutation group  is a proper subgroup of the
monodromy group of $f$. 

The fact that the removal of one extra point can change the monodromy group, makes all proofs more complicated.

{\it Thus even simplest elementary functions can have dense singular
sets and monodromy groups of cardinality of the continuum. In addition the removal of an extra point can change their monodromy groups.}

\subhead {6. Class of $\Cal S$-functions}
\endsubhead
In this section, we define a broad class of functions of one
complex variable needed in the construction of  topological  Galois theory.

\definition {Definition} A multivalued analytic function of one complex
variable is called a
$\Cal S$-function, if the set
of its singular points is at most countable. 
\enddefinition
Let us make this definition more precise. Two regular germs $f_a$ and $g_b$
defined at points $a$ and $b$ of the Riemann sphere $\Bbb S^2$ are
called {\it equivalent}
 if the germ $g_b$ is obtained from the
germ $f_a$ by the analytic continuation along some
path. Each germ $g_b$ equivalent to the germ $f_a$ is also
called a regular germ of the multivalued analytic function $f$
generated by the germ $f_a$.

A point $b\in \Bbb S^2$ is said to be a
{\it singular point}
 for the germ $f_a$ if there exists a path
$\gamma: [0,1]\to \Bbb S^2$, $\gamma (0)=a$, $\gamma (1)=b$ such
that the germ has no analytic continuation along this path, but
for any $\tau$, $0\leq \tau <1$, it admits an analytic continuation
along the truncated path $\gamma:[0,\tau]\rightarrow S^2$. 

It is easy to see that equivalent germs have the same set of singular
points. A regular germ is called a {\it $\Cal S$-germ}, if
the set of its singular points is at most countable. A
multivalued analytic function is called a $\Cal S$-function if each its regular germ is a $\Cal S$-germ.

\proclaim{Theorem 3 (on stability of the class of $\Cal S$-functions)}
  The class $\Cal S$ of all $\Cal S$-functions is stable under the following
operations:

1) differentiation, i.~e. if $f\in \Cal S$, then $f'\in
    \Cal S$;

2) integration, i.~e. if $f\in \Cal S$ and $g'=f$,
    then $g\in \Cal S$;

3) composition, i.~e. if $g,\ f\in \Cal S$, then
    $g\circ f\in \Cal S$;

4) meromorphic operations, i.~e. if $f_i\in \Cal S$,
    $i=1$, $\dots$, $n$, the function $F(x_1, \dots,x_n)$
    is a meromorphic function of $n$ variables, and
    $f=F(f_1,\dots ,f_n)$, then $f\in \Cal S$;

5) solving algebraic equations, i.~e. if $f_i\in
    \Cal S,\ i=1, \dots ,n$, and $f^n+f_1f^{n-1}+\dots
    +f_n=0$, then $f\in \Cal S$;

6) solving linear differential equations, i.e. if
    $f_i\in \Cal S,\ i=1,\dots ,n$, and
    $f^{(n)}+f_1f^{(n-1)}+\dots +f_nf=0$, then $f\in
    \Cal S$.
\endproclaim

\remark{ Remark} Arithmetic operations and the
exponentiation are examples of meromorphic operations, hence
{\it the class of $\Cal S$-functions is stable under the
arithmetic operations and the exponentiation.}
\endremark

\proclaim {Corollary 4 (see [2])}
  If a multivalued function $f$ can be obtained from single
valued $\Cal S$-functions by integration, differentiation,
meromorphic operations, compositions, solutions of algebraic
equations and linear differential equations, then the function
$f$ has at most countable number of singular points. 
\endproclaim
\proclaim {Corollary 5} A function having uncountably many singular points
cannot be represent  by generalized quadratures. In particular it cannot be a generalized  elementary function and it cannot be represented by $k$-quadratures or by quadratures.
\endproclaim

\example{Example} Consider a discrete group $\Gamma$ of fractional linear transformations  of the open unit ball $U$ having a compact fundamental domain. Let $f$ be a nonconstant  meromorphic function on $U$ invariant under the action of $\Gamma$. Each point on the boundary $
\partial U$ belongs to the closure of the set of poles  of $f$, thus the set  $\Sigma$ of singular points of $f$ contains $\partial U$. So $\Sigma$ has the cardinality of the continuum and $f$ cannot be expressed by generalized quadratures.
\endexample

\subhead {7. Monodromy group of a $\Cal S$-function}
\endsubhead
The  {\it  monodromy group } of a $\Cal{S}$-function $f$  is the group of all  permutations of
the sheets of the Riemann surface of $f$  which are induced by
motions  around the singular set $\Sigma $ of the function $f$. Below we discuss this definition more precisely. 

Let $F_{x_0}$ be the set of all  germs of the
$\Cal S$-function $f$ at point $x_0\notin \Sigma$. Consider a closed curve  $\gamma $ in
$\Bbb S^2\setminus \Sigma $ beginning and ending at the point $x_0$. Given a germ $y\in F_{x_0}$ we can continue it along the loop $\gamma$ to obtain another germ $y_\gamma\in Y_{x_0}$.
Thus each such loop $\gamma$   corresponds to a permutation $S_{\gamma}:F_{x_0}\rightarrow F_{x_0}$ of the set $F_{x_0}$ that maps a germ $y\in F_{x_0}$ to the germ $y_{\gamma}\in F_{x_0}$.

It is  easy to see that the map $\gamma\rightarrow S_\gamma$ defines a homomorphism   from the fundamental
group $\pi _1(\Bbb S^2\setminus \Sigma,x_0)$ of the domain $\Bbb S^2\setminus \Sigma $ with the base point $x_0$ to the group $S(F_{x_0})$ of permutations. The {\it
monodromy group} of the $\Cal S$-function $f$  is the image of
the fundamental group   in the group
$S(F_{x_0})$ under this homomorphism.

\remark{Remark} Instead of the point $x_0$ one can choose any other point $x_1\in \Bbb S^2\setminus \Sigma$. Such a choice will not change the    monodromy group  up to an isomorphism. To fix this isomorphism one can choose any curve $\gamma:I\rightarrow \Bbb C^N\setminus \Sigma$ where $I$ is the segment $ 0\leq t\leq 1$ and $\gamma(0)=x_0$, $\gamma(1)=x_1$ and identify each germ $f_{x_0}$ of $f$   with its continuation $f_{x_1}$ along $\gamma$.

\subhead {8. Strong non  representability  by quadratures}\endsubhead
One can prove the following  theorem.
\proclaim {Theorem 6 (see [2])}  The class of all
$\Cal S$-functions, having a solvable monodromy group, is stable 
under composition, meromorphic
operations,  integration and   differentiation.
\endproclaim
\definition{Definition} A function $f$ is {\it strongly non representable by quadratures }  if it does not belong to a class of functions defined by the following data. List of basic functions: basic elementary functions and all single valued $\Cal S$-function.
List of admissible operations: compositions, meromorphic operations, differentiation and integration.

\enddefinition

Theorem 6 implies  the  following corollary.

\proclaim {Result on  quadratures} If the monodromy group of an $S$-function $f$ is not solvable, then $f$ is   strongly non representable by quadratures.
\endproclaim

\example {Example} The monodromy group of an algebraic function $y(x)$ defined by an equation $y^5+y -x=0$ is the unsolvable group $S_5$. Thus $y(x)$ provides an example of a function with  finite set of singular points, which is strongly non representable by quadratures.
\endexample

The following corollary 7  contains  a  stronger result on non representability of algebraic functions by quadratures.
\proclaim {Corollary 7} If an algebraic function of one complex variable  has unsolvable monodromy group then it is strongly non representable by quadratures.

\endproclaim

For algebraic functions of several complex variables there is a result  similar   to Corollary 7.

\subhead {9. The  monodromy pair}
\endsubhead
The  monodromy group of a function $f$ is not only an abstract group
but is also a transitive group of permutations  of  germs of $f$ at a non singular point $x_0$.  

\definition {Definition} The {\it monodromy  pair of an $\Cal S$-function $f$} is a
pair of groups, consisting of the monodromy group  of $f$ at  $x_0$ 
 and the stationary subgroup  of a certain germ   of $f$ at $x_0$.
\enddefinition

The  monodromy pair is well defined, i.e. this pair
of groups, up to  isomorphisms, does not  depend on the  choice of the non singular point  and on the  choice of the germ of $f$ at this point.

 \definition {Definition} A pair of groups $[\Gamma,\Gamma_0]$ is an {\it almost normal pair} if there is a normal subgroup $H$ of $\Gamma$ such that $H\subset \Gamma_0$ and the coset $\Gamma_0/H$ is finite.
 \enddefinition

\definition {Definition} The pair of groups $[\Gamma,\Gamma_0]$ is called
an {\it almost solvable pair of groups} if there exists a sequence
of subgroups $$ \Gamma=\Gamma_1\supseteq\dots \supseteq \Gamma_m,
\quad \Gamma_m \subset \Gamma_0, $$ such that  for every   $i,
1\leq i\leq m-1$ group $\Gamma _{i+1}$ is a  normal divisor of
group  $\Gamma _i$ and the factor group  $\Gamma _i/\Gamma
_{i+1}$ is either a commutative group, or a finite group.
\enddefinition

\definition {Definition} The pair of groups $[\Gamma,\Gamma_0]$ is called
a {\it $k$-solvable pair of groups} if there exists a sequence
of subgroups $$ \Gamma=\Gamma_1\supseteq\dots \supseteq \Gamma_m,
\quad \Gamma_m \subset \Gamma_0, $$ such that  for every   $i,
1\leq i\leq m-1$ group $\Gamma _{i+1}$ is a  normal divisor of
group  $\Gamma _i$ and the factor group $\Gamma _i/\Gamma
_{i+1}$ is either a commutative group, or a subgroup of the group $S_k$ of permutations of $k$ elements.
\enddefinition

We say that  group $\Gamma $ is
{\it almost solvable} or {\it $k$-solvable} if pair $[\Gamma ,e]$, where $e$ is the  group
containing only the unit element, is almost solvable or $k$-solvable respectively. 

It is easy to see that {\it an almost normal pair of groups $[\Gamma,\Gamma_0]$ is almost solvable or $k$-solvable if and only if the group $\Gamma$ is almost solvable or $k$-solvable respectively}.

\subhead {10. Strong non  representability  by $k$-quadratures}\endsubhead
One can prove the following  theorem.
\proclaim {Theorem 8 (see [2])}  The class of all
$\Cal S$-functions, having a $k$-solvable monodromy pair, is stable
under composition, meromorphic
operations,  integration, differentiation and solutions of algebraic equations of degree $\leq k$.
\endproclaim
\definition{Definition}  A function $f$ is {\it strongly  non representable by $k$-quadratures}  if it does not belong to a class of functions defined by the following data. List of basic functions: basic elementary functions and all single valued $\Cal S$-function.
List of admissible operations: compositions, meromorphic operations, differentiation,integration and solutions of algebraic equations of degree $\leq k$.
\enddefinition 
Theorem 8 implies  the  following corollary.

\proclaim {Result on  $k$-quadratures} If the monodromy pair of an $S$-function $f$ is not $k$-solvable, then $f$ is strongly  non representable by $k$-quadratures.
\endproclaim

\example {Example} The monodromy group of an algebraic function $y(x)$ defined by an equation $y^n+y -x=0$ is the permutation group group $S_n$. For $n\geq 5$ the group $S_n$ is not an $(n-1)$-solvable group. Thus $y(x)$ provides an example of a function with  finite set of singular points which is strongly non representable  by $(n-1)$-quadratures.
\endexample

This example can be  generalized.
\proclaim {Corollary 9 (see [2])} If an algebraic function of one complex variable  has non $k$-solvable monodromy group then it is strongly non representable by $k$-quadratures.

\endproclaim

\proclaim {Theorem 10 (see [2])} An  algebraic function of one variable whose monodromy group is $k$-solvable, can be represented by $k$-radicals.
\endproclaim
Results similar to Corollary 9 and Theorem 10 hold also for algebraic functions of several complex variables.

\subhead {11. Strong non  representability  by generalized quadratures}\endsubhead
One can prove the following  theorem.
\proclaim {Theorem 11 (see [2])}   The class of all
$\Cal S$-functions, having an almost solvable monodromy pair, is stable
under composition, meromorphic
operations,  integration, differentiation and solutions of algebraic equations.
\endproclaim

\definition{Definition}  A function $f$ is {\it strongly  non representable by generalized quadratures}  if it does not belong to a class of functions defined by the following data. List of basic functions: basic elementary functions and all single valued $\Cal S$-function.
List of admissible operations: compositions, meromorphic operations, differentiation, integration and solutions of algebraic equations.
\enddefinition
Theorem 11 implies  the  following corollary.

\proclaim {Result on  generalized quadratures} If the monodromy pair of an $S$-function $f$ is not almost solvable, then $f$ is strongly  non representable by generalized quadratures.
\endproclaim

Suppose that the Riemann surface of a
function $f$ is a universal covering space over the Riemann sphere with $n$ punched points. If $n\geq 3$ then the  function $f$ is strongly non representable
by generalized quadratures. Indeed,
the monodromy pair of  $f$ consists of the free  group with $n-1$ generators, and its unit subgroup. It is easy to see that
such a pair of groups is not almost solvable.

\example{Example} Consider the function $z(x)$, which maps  the upper half-plane onto a triangle with vanishing angles, bounded by three circular arcs.
The Riemann surface of  $z(x)$ is a universal covering space over the
sphere with three punched points. \footnote {it is easy to see that the function $z(x)$ maps its Riemann surface to the open ball whose boundary contains the vertices of the triangle. These properties of the function $z(x)$ play the crucial role in  Picard's beautiful proof of his Little Picard Theorem.}  Thus $z(x)$ is strongly non
representable by generalized quadratures.
\endexample

\example{Example}
Let $K_1$ and $K_2$ be the following elliptic
integrals, considered as the functions of the parameter $x$:
  $$
K_1(x)=\int_0^1\frac{dt}{\sqrt{(1-t^2)(1-t^2x^2)}}\ {\hbox{ and }}
K_2(x)=\int_0^{\frac{1}{x}}\frac{dx}{\sqrt{(1-t^2)(1-t^2x^2)}}. $$
The functions  $z(x)$ can be obtained
from $K_1(x)$ and from $K_2(x)$  by quadratures. Thus both functions
 $K_1(x)$ and $K_2(x)$ are strongly non representable  by generalized quadratures.
\endexample

In the next section  we will list all polygons $G$ bounded by circular arcs for which the Riemann map of the upper half-plan onto $G$ is representable by generalized quadratures.

\subhead{12. Maps of the upper half-plane onto a curved polygon} 
\endsubhead
Consider a  polygons $G$ on the complex plane  bounded by circle arcs, and the function $f_G$ establishing the Riemann mapping of the upper half-plane onto the polygon
$G$. The Riemann--Schwarz reflection
principle allows to describe the monodromy group $L_G$ of the function $f_G$ and to show that all singularities of $f_G$ are simple enough. This information together with Theorem  11 provide a complete classification of all polygons $G$ for which the function $f_G$ is representable in explicit form (see~[2]).

If a polygon $\tilde G$ is obtained from a polygon $G$ by a ]linear transformation $w:\Bbb C  \to \Bbb C$ then $f_{\tilde G}=w (f_G)$. Thus it is enough to classify $G$ up to a ]linear transformation.

{\bf 1) The first case of integrability:}  the continuations of all sides of the polygon $G$ intersect at one point.

Mapping this point to infinity by a fractional linear transformation, we obtain a polygon $G$ bounded by straight line segments.
All transformations in the group $L(G)$ have the form $z \to az+b.$ All germs of
the function $f = f_G$ at a non-singular point c are obtained from a fixed germ $f_c$ by
the action of the group $L(G)$ consisting of the affine transformations $z\to a z+b$. The germ $R_c = (f''/ f)_c$ is invariant under the action of the group $L(G)$. Therefore, the germ $R_c$ is a germ of a single valued function $R$.
The singular points  of  $R$ can only be poles (see ). Hence the function $R$ is  rational. The equation $ f''/ f= R$ is
integrable by quadratures. This integrability case is well known. The function $f $ in
this case is called the {\it Christoffel--Schwarz integral}.

{\bf 2) The second case of integrability:} there is  a pair of points such that, for every side of the polygon
G, these points are either symmetric with respect to this side or belong to the continuation
of the side.

We can map these two points to zero and infinity by a fractional
linear transformation. We obtain a polygon $G$ bounded by circle arcs centered at
point 0 and intervals of straight rays emanating from 0 (see Figure 2). All transformations
in the group $L(G)$ have the form $z\to az, z \to b/z.$  All germs of the function
$f = f_G$ at a non-singular point $c$ are obtained from a fixed germ $f_ c $ by the action of
the group $L(G)$ :$$f_ c \to a  f_c, f_ c \to b/ f _c.$$ The germ $R_c = (f'_ c/ f_ c)^2$ is invariant under the action of the group $L(G)$. Therefore, the germ $R_c$ is a germ of a single valued function $R$. The singular points  of  $R$ can only be poles (see ). Hence the function $R$ is  rational. The equation $R = (f'/ f)^2$ is integrable by quadratures
\bigskip

{\bf 3) The finite nets of  circles}. To describe the third  case of integrability we need to define first the finite net of circles on the complex plane. The classification of finite groups, generated by reflections in the Euclidian space $\Bbb R^3$ is well known. Each such group is the symmetry  group of the following bodies:

1. a regular n-gonal pyramid;

2. a regular n-gonal diheron, or the body formed by two equal regular n-gonal pyramids
sharing the base;

3. a regular tetrahedron;

4. a regular cube or icosahedron;

5. a regular dodecahedron or icosahedron.

All these groups of isometries, except for the group of dodecahedron or icosahedron,
are solvable.

The intersections of the unit sphere, whose center coincides with the
barycenter of the body, with the mirrors, in which the body is symmetric, is a certain
net of great circles.  Stereographic projections of each of them is a net net of circles on complex plane defined up to a fractional linear transformation.  The nets corresponding to the bodies listed above will be called
the finite nets of  circles.
\bigskip
{\bf 4) The third case of integrability:} every side side of a polygon $G$  belongs to some finite net  of  circles. In this case  the function $f_G$ has finitely many branches. Since all singularities of
the function $f_G$ are algebraic (see [2]), the
function $f_G$ is an algebraic function. For all finite nets but the net  of dodecahedron or icosahedron, the algebraic function $f_G$ is representable by radicals. For the net  of dodecahedron or icosahedron the function $f_G$ is representable by radicals and solutions of degree five algebraic equations (in other words $f_G$ is representable by $k$-radicals).
\medskip

{\bf 5) The strong non representability}. Our results  imply the following:
\proclaim {Theorem 12 (see [2])} If a  polygon $G$ bounded by circles arcs  does not belong to one of
the three  cases described  above, then the function $f_G $ is strongly non  representable by generalized quadratures.
\endproclaim

\subhead {13. Non solvability of linear differential equations}\endsubhead 
Consider a homogeneous linear differential equation
$$ 
y^{(n)}+ r_1y^{(n-1)}+\dots +r_ny=0,\tag 3
$$
whose coefficients  $r_i$'s are rational functions  of the complex variable
$x$.  The set  $\Sigma\subset \Bbb C$  of poles of  $r_i$'s is called  {\it the set of
singular points}   of the equation (3). At a  point $x_0\in \Bbb C\setminus \Sigma$ the germs of solutions of (3) form a $\Bbb C$-linear  space $V_{x_0}$ of dimension $n$. The {\it monodromy group $M$ of the equation (3)} is the group of all linear transformations  of the space $V_{x_0}$  which are induced by motions around the set $\Sigma$. Below we discuss this
definition more precisely.

Consider a closed curve  $\gamma $ in
$\Bbb C\setminus \Sigma $ beginning and ending at the point $x_0$. Given a germ $y\in V_{x_0}$ we can continue it along the loop $\gamma$ to obtain another germ $y_\gamma\in V_{x_0}$.
Thus each such loop $\gamma$   corresponds to a map $M_{\gamma}:V_{x_0}\rightarrow V_{x_0}$
 of the space $V_{x_0}$to itself  that maps a germ $y\in V_{x_0}$ to the germ $y_{\gamma}\in V_{x_0}$. The map $M_\gamma$ is  linear  since an analytic continuation respects the arithmetic operations.
It is easy to see that the map $\gamma\rightarrow M_\gamma$ defines a homomorphism   of the fundamental
group $\pi _1(\Bbb C\setminus \Sigma,x_0)$ of the domain $\Bbb C\setminus \Sigma $ with the base point $x_0$ to the group $GL(n)$ of invertible linear transformations  of  the space  $V_{x_0}$. 

The {\it
monodromy group} $M$ of the equation (3) is the image of
the fundamental group   in the group
$GL(n))$ under this homomorphism.

\remark{Remark} Instead of the point $x_0$ one can choose any other point $x_1\in \Bbb C\setminus \Sigma$. Such a choice will not change the    monodromy group  up to an isomorphism. To fix this isomorphism one can choose any curve $\gamma:I\rightarrow \Bbb C^N\setminus \Sigma$ where $I$ is the segment $ 0\leq t\leq 1$ and $\gamma(0)=x_0$, $\gamma(1)=x_1$ and identify each germ $y_{x_0}$ of solution of (*) with its continuation $y_{x_1}$ along $\gamma$.
\endremark

\proclaim {Lemma 13}  The stationary subgroup in the monodromy group $M$ of the germ $y\in V_{x_0}$   of almost every solution  of the equation (3)  is trivial (i.e. contains only the unite  element $e\in M$).
\endproclaim
\demo {Proof} The monodromy group $M$ contains countable many linear transformations $M_i$. The space $L_i\subset V_{x_0}$ of fixed points of a non identity transformation $M_i$,  is a proper subspace of $ V_{x_0}$. The union $L$ of all subspaces $L_i$ is a measure zero subset of  $V_{x_0}$. The stationary subgroup in  $M$ of $y\in V_{x_0}\setminus L$  is trivial.
\enddemo

\proclaim {Theorem 14 (see [2])}  If the monodromy group of the equation   (3)
is not almost solvable (is not $k$-solvable, or is not solvable)  then its almost every  solution is strongly  non representable by generalized quadratures  (correspondingly, is strongly non representable by $k$-quadratures, or is strongly non representable by quadratures).

\endproclaim

Consider a homogeneous system of linear differential equations  

$$
 y'=A y \tag 4
$$
where $y=(y_1,\dots,y_n)$ is the unknown vector valued function and $A=\{a_{i,j}(x)\}$
is a $n\times n$ matrix, whose entries are rational functions of the
complex variable $x$. One can define the {\it monodromy group} of the equation (4) exactly in the same way as it was defined for the equation (3).

We will say that a vector valued function $y=(y_1,\dots, y_n)$  belongs to a certain class of functions if all its components $y_i$  belong to this class. For example the statement "a vector valued function $y=(y_1,\dots,y_n)$ is strongly non representable by generalized quadratures" means that at least one component $y_i$ of $y$ is strongly non representable by generalized quadratures.

\proclaim {Theorem 15}  If the monodromy group of the system  (4)
is not almost solvable (is not $k$-solvable, or is not solvable)  then its almost every  solution is strongly  non representable by generalized quadratures  (correspondingly, is strongly non representable by $k$-quadratures, or is strongly non representable by quadratures).

\endproclaim

\subhead{14. Solvability of  Fuchsian  equations}
\endsubhead
The differential field  of rational functions of $x$ is isomorphic to the differential field $\Cal R$ of germs of rational functions at the point $x_0\in \Bbb C\setminus \Sigma $.
Consider the differential field extension  $\Cal R \{y_1,\dots,y_n\}$ of $\Cal R$ where the  germs  $ y_1,\dots y_n$ form a basis in the space $V_{x_0}$ of solutions of the equation~(3) at $x_0$.

\proclaim {Lemma 16} Every linear map $M_\gamma$ from the monodromy group of equation (3), can be uniquely extended to a differential automorphism of the  differential field $\Cal R\{ y_1,\dots,y_n\}$ over the field $\Cal R$.
\endproclaim

\demo {Proof} Every element $f\in \Cal R \{y_1,\dots,y_n\}$ is a rational function of the independent variable $x$, the germs of solutions $y_1,\dots,y_n$ and their derivatives. It can be continued meromorphically along  the curve
$\gamma\in \pi_1(\Bbb C\setminus \Sigma,x_0)$  together with  $y_1,\dots,y_n$.
This continuation gives the required differential
automorphism, since the continuation preserves the arithmetical
operations and differentiation, and every rational function of $x$ returns back
to its original  values (since it is a single-valued valued function). The differential automorphism is unique because the extension is generated by $y_1,\dots, y_n$.
\enddemo

The {\it differential Galois group} (see [2], [3]) of the equation (3) over $\Cal R$ is the group of all differential automorphisms  of the  differential field $\Cal R\{ y_1,\dots, y_n\}$ over the differential field  $\Cal R$. According to Lemma 32 the monodromy group  of the equation (3) can be considered as a subgroup of its differential Galois group  over $\Cal R$.

The differential  field of  invariants of the monodromy group action is a subfield
of $\Cal R\{y_1,\dots,y_n\}$, consisting of the single-valued
functions. Differently from the algebraic case, for differential
equations  the field of invariants under the action of the
monodromy group can be bigger than the field of rational
functions. The reason is that  the solutions
of differential equations may grow exponentially in approaching
the singular points or infinity.

\example {Example} All  solutions of the simplest  differential equation $y'=y$  are
single-valued exponential functions $y=C\exp x$, which are not rational.
\endexample

For a wide class of   Fuchsian  linear
differential equations all the solutions,  while approaching
the singular points and the point infinity, grow polynomially.

The following  Frobenius theorem is an analog for Fuchsian  equations of  C.Jordan theorem (see [ ]) for algebraic equations.

\proclaim {Theorem (Frobenius)}  For Fuchsian differential  equations 
 the subfield of the differential field
$\Cal R\{y_1,\dots,y_n\}$, consisting of single-valued
functions, coincides with the field of rational functions.
\endproclaim

A system of linear differential equations (4) is called a {\it Fuchsian system} if the matrix  $A$ has the  following form:
 $$
 A(x)=\sum _{i=1}^k\dfrac{A_i}{x-a_i},\tag 5
 $$

where the $A_i$'s are constant matrices. Linear Fuchsian system of differential equations in many ways are similar to  linear Fuchsian  differential equations

In construction of explicit solutions of linear  differential equations the following theorem is needed.

\proclaim{Theorem (Lie--Kolchin)} Any  connected  solvable algebraic group acting by linear transformations on  a finite-dimensional vector space over  $\Bbb C$ is  triangularizable in a suitable basis.
\endproclaim 

Using Frobenius Therem and Lie--Kolchin Theorem  one can prove that the only
reasons for unsolvability of  Fuchsian linear
differential equations and systems of linear differential equations are topological. In other words, if there are no topological
obstructions to solvability then such equations and systems of equations are solvable. Indeed, the following theorems hold:

\proclaim {Theorem 17 (see [2])}  If the monodromy group of the linear Fuchsian differential  equation (3)
is almost solvable (is  $k$-solvable, or is solvable)  then  its every  solution is  representable by generalized quadratures  (correspondingly, is  representable by $k$-quadratures, or is  representable by quadratures).

\endproclaim

\proclaim {Theorem 18 (see [2])}  If the monodromy group of the linear Fuchsian system differential equations (4)
is  almost solvable (is  $k$-solvable, or is  solvable)  then  its every  solution is    representable by generalized quadratures  (correspondingly, is  representable by $k$-quadratures, or is  representable by quadratures).

\endproclaim

\subhead{15.   Fuchsian systems  
with small coefficients}
\endsubhead
In general the monodromy group of a given Fuchsian equation is very hard to compute.  It is known only for very special equations, including the famous hypergemetric equations. Thus Theorems 17 and 18 are not explicat. 

If the matrix $A(x)$ in the system (4) is triangular then one can easily  solve the system by quadratures. It turns out that if the matrix $A(x)$ has the form (5), where the matrices $A_i$'s are sufficiently small, then the system (4) with a non triangular matrix $A(x)$ is unsolvable by generalized quadratures for a topological reason. 

\proclaim {Theorem 19 (see [2])} If the  matrices $A_i$'s are sufficiently small, $\Vert
 A_i\Vert<\varepsilon(a_1,\dots,a_k,n)$, then the  monodromy group of the system
 $$ y'= (\sum _{i=1}^k\dfrac{A_i}{x-a_i})y \tag 6$$
is almost solvable if and only if  the  matrices $A_i$'s are triangularizable in a suitable basis.
\endproclaim
\proclaim {Corollary 20} If in the assumptions of Theorem 19 the  matrices $A_i$'s are not triangularizable in a suitable basis then almost every  solution of the system (6) is strongly non   representable by generalized quadratures.
\endproclaim

\subhead {16. Polynomials invertible by radicals}
\endsubhead
In 1922 J.F.Ritt published (see [5]) the following beautiful theorem which fits nicely into topological Galois theory.
\proclaim {Theorem (J.F. Ritt)}
The inverse function of a polynomial with complex coefficients can be represented by radicals if and only if the polynomial is a composition of linear polynomials, the power polynomials $z\to z^n$, Chebyshev polynomials and polynomials of degree at most 4.
\endproclaim

\demo {Outline of proof  (following  [6])}

1)   {\it Every polynomial is a composition of primitive ones:} Every polynomial is a composition of polynomials that are not themselves compositions of polynomials of degree $>1$. Such polynomials are called {\it primitive}. Recall that a permutation group $G$ acting on a non-empty set $X$ is called {\it primitive} if $G$ acts transitively on $X$ and G preserves no nontrivial partition of $X$. {\it A polynomial is primitive if and only if the monodromy group of inverse of the polynomial   acts primitively on its branches}.
\smallskip

2)  {\it Reduction to the case of primitive polynomials:} A composition of polynomials is invertible by radicals if and only if each polynomial in the composition is invertible by radicals. Indeed, if each of the polynomials in composition is invertible by radicals, then their composition also is. Conversely, if a polynomial $R$ appears in the presentation of a polynomial $P$ as a composition $P=Q\circ R\circ S$ and $P^{-1}$ is representable by radicals, then $R^{-1}=Q\circ P^{-1}\circ S$ is also representable by radicals. Thus it is enough to classify only the primitive polynomials invertible by radicals.
\smallskip

3) {\it A result on solvable primitive permutation groups containing a full cycle:} A primitive polynomial is invertible by radicals if and only if the monodromy group of inverse of the polynomial is solvable.  Since it acts primitively on its branches  and contains a full cycle (corresponding to a loop around the point at infinity on the Riemann sphere), the following group-theoretical result of Ritt is useful for the classification of polynomials invertible by radicals:
\proclaim {Theorem (on primitive solvable groups with a cycle)}
Let $G$ be a primitive solvable group of permutations of a finite set $X$ which contains a full cycle. Then either $|X|=4$, or $|X|$ is a prime number $p$ and $X$ can be identified with the elements of the field $F_p$ so that the action of $G$ gets identified with the action of the subgroup of the affine group $AGL_1(p)=\{x\to ax+b|a\in (F_p)^*,b\in F_p\}$ that contains all the shifts $x\to x+b$.
\endproclaim

4) {\it Solvable monodromy groups of inverse of primitive polynomials:} It can be shown by applying Riemann--Hurwitz formula that among the groups in Theorem on primitive solvable groups with a cycle,  only the following groups can be realized as monodromy groups of inverse of primitive polynomials: 1. $G\subset S_4$, 2. Cyclic group $G=\{x\to x+b\}\subset AGL_1(p)$, 3. Dihedral group $G=\{x\to \pm x+b\}\subset AGL_1(p)$.
\smallskip

5) {\it Description of primitive polynomials invertible by radicals:} It can be easily shown (see for instance [see Ritt 22], [Khovanskii 07 Variations], [Burda Khovanskii11 Branching]) that the following result holds:
\proclaim{Theorem 21}
If the monodromy group of  inverse of a primitive polynomial  is a subgroup of the group $\{x\to \pm x+b\}\subset AGL_1(p)$, then up to a linear change of variables the polynomial is either a power polynomial or a Chebyshev polynomial.
\endproclaim

Thus the polynomials whose inverse have monodromy groups 1-3 are respectively 1. Polynomials of degree four. 2. Power polynomials up to a linear change of variables. 3. Chebyshev polynomials up to a linear change of variables.

In each of these cases the fact that the polynomial is invertible by radicals follows from solvability of the corresponding monodromy group or from explicit formulas for its inverse (see for instance [BurdaKhovanskii11Branching]).
\enddemo
\subhead {17. Polynomials invertible by $k$-radicals}
\endsubhead
In this section we discuss the following generalization of J.F.Ritt's Theorem.  

\proclaim {Theorem 22 (see [6])}
 A polynomial invertible by radicals and solutions of equations of degree at most $k$ is a composition of power polynomials, Chebyshev polynomials, polynomials of degree at most $k$ and, if $k\leq 14$, certain primitive polynomials whose inverse have  exceptional monodromy groups. A description of these exceptional polynomials can be given explicitly.
\endproclaim

The proofs rely on classification of monodromy groups of inverse of primitive polynomials  obtained by M\"{u}ller based on group-theoretical results of Feit and on previous work on primitive polynomials whose inverse have exceptional monodromy groups by many authors. Besides the references to these highly involved and technical results an outline of the proof of Theorem 22 is not complicated and it resembles the outline of the proof of Ritt's Theorem.

Let us start with some background on representability by $k$-radicals.

\definition{Definition}
Let $k$ be a natural number. A field extension $L/K$ is $k$-radical if there exists a tower of extensions $K=K_0\subset K_1\subset \ldots \subset K_n$ such that $L\subset K_n$ and for each $i$, $K_{i+1}$ is obtained from $K_i$ by adjoining an element $a_i$, which is either a solution of an algebraic equation of degree at most $k$ over $K_i$, or satisfies $a_i^m=b$ for some natural number $m$ and $b\in K_i$.
\enddefinition

\proclaim {Theorem 23 (see [2])} A Galois extension $L/K$ of fields of characteristic zero is $k$-radical if and only if its Galois group is $k$-solvable.
\endproclaim

An algebraic function $z=z(x)$ of one or several complex variables is said to be {\it representable by  $k$-radicals} if the corresponding extension of the field of rational functions is a $k$-radical extension.

Theorem 23 and C. Jordan's Theorem (see [ 1]) imply the following corollary.
 
 \proclaim {Corollary 24}  An algebraic function is representable by $k$-radicals if and only if its monodromy group is $k$-solvable.
 \endproclaim
 (Note that Theorem 10 above coincides with a  part of  Corollary 23).

Let us outline briefly the main steps in the proof of Theorem 22:

\demo {Outline of proof of Theorem 22}:
 
1) Exactly as  in Ritt's theorem one can show that a composition of polynomials is invertible by $k$-radicals if and only if each polynomial in the composition is invertible by $k$-radicals. Thus one can reduce Theorem 22  to the case of primitive polynomials.
\smallskip

2) Feit and Jones totally classified all primitive permutation groups of $n$ elements containing a full cycle.
\smallskip 

3) Using the this classification  and Riemann-Hurwitz formula, M\"{u}ller listed
all groups of permutations of $n$ elements  which are monodromy groups of inverses of degree $n$ primitive polynomials.
\smallskip
4) For each group from M\"{u}ller's list of groups of permutations of $n$ elements one can determine the smallest $k$ for which it is
$k$-solvable and choose the {\it exceptional groups} for which $k$ is smaller than $n$. 
\smallskip
5) For each such exceptional group one can explicitly describe polynomials whose inverse has the exceptional monodromy group.
\enddemo

\subhead{18. Acknowledgement}\endsubhead
I would like to thank Michael Singer who invited me to write comments for a new edition of the classical J.F.~Ritt's book ``Integration in finite terms''~[4]. This preprint was written as a part of these comments. I also am grateful to Fedor Kogan who edited my English.
\bigskip

\centerline {REFERENCES}
\medskip

[1] A.G. Khovanskii, On representability of algebraic functions by radicals,\linebreak arXiv:1903.08632 [math.AG]

[2]  A. Khovanskii, Topological Galois theory. Solvability and unsolvability of equations in
finite terms. Translated by Valentina Kiritchenko and Vladlen Timorin. Series: Springer
Monographs in Mathematics. Springer Berlin Heidelberg. 2014, XVIII, 305 pp. 6 illus.

[3] M. van der Put, M. Singer, Galois theory of linear differential equations (Springer, Berlin/New York, 2003).

[4] J. Ritt, Integration in finite terms. Liouville's theory of elementary methods, N. Y. Columbia Univ. Press. 1948.

[5] J. Ritt, On algebraic functions which can be expressed in terms of radicals, Trans. Am. Math. Soc. 24, 21--30 (1922).

[6] Yu.Burda, A.Khovanskii, Polynomials invertible in $k$-radicals. Arnold Mathematical Journal, V. 2, No 1, 2016, 121--138.

\end